\documentclass[a4paper,11pt]{article}
\usepackage{graphicx,multirow}
\usepackage{amssymb,amsmath,amsthm}
\usepackage{lscape,lipsum,microtype}
\usepackage[margin=2.54cm]{geometry}  
\usepackage{hyperref}

\usepackage{ifluatex}
\ifluatex
 \usepackage{unicode-math}
\fi

\usepackage{latexsym}
\usepackage{amsmath}
\usepackage{amsfonts}
\usepackage{amssymb}
\usepackage{amscd}
\usepackage{amsmath}

\newcommand{\D }{\Delta }

\renewcommand{\l }{\lambda }

\newcommand{\n }{\nabla }

\newcommand{\intbar}{\mathop{\int\makebox(-13.5,0){\rule[4pt]{.7em}{0.3pt}}%
\kern-6pt}\nolimits}

\newcommand{\be}{\begin{equation}}
\newcommand{\ee}{\end{equation}}
\newcommand{\bes}{\begin{equation*}}
\newcommand{\ees}{\end{equation*}}
\newcommand{\ba}{\begin{eqnarray}}
\newcommand{\ea}{\end{eqnarray}}
\newcommand{\bas}{\begin{eqnarray*}}
\newcommand{\eas}{\end{eqnarray*}}
\newenvironment{pf}{\noindent{\sc Proof}.\enspace}{\rule{2mm}{2mm}\medskip}
\newenvironment{pfn}{\noindent{\sc Proof}}{\rule{2mm}{2mm}\medskip}

\newcommand{\R}{\mathbb{R}}

\author{ Cheikh Birahim NDIAYE}

\date{}

\title{\bf Optimal control for the conformal Laplacian obstacle problem}
\begin{document}

\newtheorem{lem}{Lemma}[section]
\newtheorem{pro}[lem]{Proposition}
\newtheorem{thm}[lem]{Theorem}
\newtheorem{rem}[lem]{Remark}
\newtheorem{cor}[lem]{Corollary}
\newtheorem{df}[lem]{Definition}

\maketitle

\begin{center}
{\small

\noindent  Department of Mathematics Howard University \\  Annex $3$, Graduate School of Arts and Sciences, $217$\\ DC 20059 Washington, USA

}

\end{center}

\footnotetext[1]{E-mail addresses: cheikh.ndiaye@howard.edu\\
\thanks{\\ The author was partially supported by NSF grant DMS--2000164.}}

\

\

\begin{center}
{\bf Abstract}
\end{center}
We study an optimal control problem associated to the conformal Laplacian obstacle problem on closed \;$n$-dimensional Riemannian manifolds with $n\geq 3$. When the Yamabe invariant of the Riemannian manifold is positive, we show that the optimal controls are equal to their associated optimal states and show the existence of a smooth optimal control which induces a conformal metric with constant scalar curvature. For the standard sphere, we prove that the standard bubbles-namely conformal factor of metrics conformal to the standard one with constant positive scalar curvature- are the only optimal controls and hence equal to their associated optimal state.

 \begin{center}


\bigskip

\centerline{\bf AMS subject classification: 53C21, 35C60, 58J60, 55N10.}

\end{center}


\section{Introduction and statement of the results}\label{intro}
The resolution of the Yamabe problem by Yamabe\cite{yamabe}, Trudinger\cite{tru}, Aubin\cite{aubin1}, and Schoen\cite{schoen} (chronologically) has been a landmark in Geometric Analysis, see also  \cite{bah1}, \cite{bahbre} for another approach. The Yamabe problem is the geometric question which ask whether every closed Riemannian manifold of dimension greater than $2$ carries a conformal metric with constant scalar curvature. Analytically it is equivalent to finding a smooth positive solution to 
\begin{equation}\label{eq:sequation}
L_gu=cu^{\frac{n+2}{n-2}}\;\;\text{on}\;\;M,
\end{equation}
where $c$ is a real number, $(M, g)$ is the background closed Riemannian manifold of dimension $n\geq3$,  and $$L_g:=-4\frac{n-1}{n-2}\D_g+R_g$$ is the conformal Laplacian with $R_g$ the scalar curvature of \;$(M, g)$ \;and \;$\D_g$ the Laplace-Beltrami associated to \;$(M, g)$.  We recall that \;$\frac{n+2}{n-2}=2^*-1$\; where \;$2^*=\frac{2n}{n-2}$.
\vspace{5pt}

\noindent
The solution of \eqref{eq:sequation} given by the works of Yamabe\cite{yamabe}, Trudinger\cite{tru}, Aubin\cite{aubin1}, and Schoen\cite{schoen} is obtained by finding a smooth minimizer of the Yamabe functional $J^g$ defined by
$$
J^g(u):=\frac{\left<u, u\right>_g}{||u||_{L^{2^*}(M, g)}}, \;\;\;\;u\in H^1_+(M, g):=\{u\in H^1(M, g), \;u>0\},
$$
where 
$$
\left<u,v\right>_{g}=4\frac{n-1}{n-2}\int_M\n_gu\n_gvdV_g+\int_{M}R_guvdV_g, \;\;u, v\in H^1_+(M, g),
$$
$dV_g$ is the Riemannian measure associated to \;$(M, g)$, $L^{2^{*}}(M, g)$ is the standard Lebesgue space of functions which are \;$2^*$-integrable over \;$M$\; with respect to \;$g$, $||\cdot||_{L^{2^*}(M, g)}$\; is the standard  \;$L^{2^{*}}$-norm on $L^{2^{*}}(M, g)$\; and \;$H^1(M, g)$\; is the space of functions on \;$M$\; which are of class \;$L^2$\; together with their first derivatives with respect to $g$, see \cite{aubin}, \cite{gt},\cite{lp} for precise definitions.
We remark that for \;$u$\; smooth,
 $$
 \left<u,u\right>_{g}=\left<L_gu, u\right>_{L^2(M, g)},
 $$
where \;$\left<\cdot, \cdot\right>_{L^2(M, g)}$\;denotes the \;$L^2$\; scalar product on \;$M$\; with respect to \;$g$.
\vspace{8pt}

\noindent
The Yamabe  problem \eqref{eq:sequation} and the prescribed scalar curvature problem ($K$ a smooth function) 
\begin{equation}\label{eq:sequationk}
L_gu=Ku^{\frac{n+2}{n-2}}\;\;\text{on}\;\;M,
\end{equation}
has been intensively studied in the susing methods of Calculus of Variations, Critical Points Theory, Morse Theory, Dynamical Systems, Blow-up Analysis, Perturbations Methods, Algebraic Topology,  see \cite{aubin1}, \cite{bc}, \cite{bcch}, \cite{bren2}, \cite{schoen}, \cite{lp},  \cite{tru}, \cite{yamabe} and the references therein.
\vspace{10pt}

\noindent
In this paper, we investigate equation \eqref{eq:sequation} in the context of Optimal Control Theory. We recall that in the works of Yamabe\cite{yamabe}, Trudinger\cite{tru}, Aubin\cite{aubin1}, and Schoen\cite{schoen}, the most difficult case is the positive one, namely \;$\mathcal{Y}(M, [g])>0$, where
$$\mathcal{Y}(M, [g]):=\inf_{u\in H^1_+(M, g) }J^g(u)$$ is the Yamabe invariant of $(M, [g])$ and $$[g]=\{g_u=u^{\frac{4}{n-2}}g,\;\; u\in C^{\infty}_+(M, g)\}$$  denotes the conformal class of $g$\; with
$$C^{\infty}_+(M, g)=\{u\in C^{\infty}(M, g):\;\;u>0\}$$ and \;$C^{\infty}(M, g)$\; is the space of \;$C^{\infty}
$-functions on $M$ with respect to $g$, see \cite{aubin}, \cite{gt}, \cite{lp} for a precise a definition.
\vspace{7pt}

\noindent
In this paper, we focus in this case. Precisely, under the assumption \;$\mathcal{Y}(M, [g])>0$, we study the following optimal control problem for the conformal Laplacian obstacle problem
\begin{equation}\label{stateop}
\text{Finding} \;\;\;u_{min}\in H^1_+(M, g)\;\;\text{such that}\;\;\; I^g(u_{min})=\min _{u\in H^1_+(M, g)}I^g(u),
\end{equation}
where
 $$
 I^g(u)= \frac{ \left<u,u\right>_{g}}{||T_g(u)||_{L^{2^*}(M, g)}^2} ,\;\;\;\;u\in H^{1}_+(M, g)
$$
with 
 \begin{equation}\label{defoo}
 T_g(u)=\arg{\min_{v\in H^1_+(M, g),\\\;\;v\geq u}\left<v, v \right>_g}
\end{equation}
where the symbol $$\arg{\min_{v\in H^1_+(M, g),\\\;\;v\geq u}\left<v, v \right>_g}$$ denotes the unique solution to the minimization problem \;$$\min_{v\in H^1_+(M, g),\\\;\;v\geq u}\left<v, v \right>_g,$$ see Lemma \ref{obstacleq}.  
\vspace{10pt}

\noindent
We obtain the following result for the optimal control problem \eqref{stateop}.
\begin{thm}\label{positivevarmass}
Assuming that  \;$\mathcal{Y}(M,\;[g])>0$, then\\
1)\\ For $u\in H^1_+(M, g)$, 
$$
I^g(u)=\min _{v\in H^1_+(M, g)}I^g(v)\implies T_g(u)=u, \;\;u\in C^{\infty}_+(M, g),\;\;\text{and}\;\;R_{g_u}=const>0,
$$
with \;$g_u=^{\frac{4}{n-2}}g$.\\
2)\\
 There exists \;$u_{min}\in C^{\infty}_+(M, g)$\; such that
$$
I^g(u_{min})=\min _{v\in H^1_+(M, g)}I^g(v)\;\;\;{and}\;\;\;R_{g_{u_{min}}}=\mathcal{Y}(M,\;[g]),
$$
with \;$g_{u_{min}}=u_{min}^{\frac{4}{n-2}}g$.
\end{thm}
\vspace{7pt}

\noindent
 When \;$(M, g)=(\mathbb{S}^n, g_{\mathbb{S}^n})$ is the standard \;$n$-dimensional sphere of \;$\R^{n+1}$, we have.
\begin{thm}\label{sphere}
Assuming that  \;$(M, g)=(\mathbb{S}^n, g_{\mathbb{S}^n})$, then for $u\in H^1_+(M, g)$,
$$
I^g(u)=\min _{v\in H^1_+(M, g)}I^g(v)\;\;\;\;\text{is equivalent to }\;\;\; u\in C^{\infty}_+(M, g)\;\;\text{and}\;\;R_{g_u}=const>0,
$$
with \;$g_u=^{\frac{4}{n-2}}g$\; and \;$R_{g_{u}}$\; the scalar curvature of \;$(M, g_u)$.
\end{thm}

\begin{rem}\label{standardbubble}
\begin{itemize}
\item
In Theorem \ref{positivevarmass}, we have point 1) $\implies$\;\;$u_{min}= T_g(u_{min})$.
\item
We recall that the standard bubbles of the Yamabe problem are the functions $u\in C^{\infty}_+(\mathbb{S}^n, g_{\mathbb{S}^n})$ with $R_{(g_{\mathbb{S}^n})_u}=const>0$\;where\;$(g_{\mathbb{S}^n})_u=u^{\frac{4}{n-2}}g_{\mathbb{S}^n}$. Thus, by point 1) of Theorem \ref{positivevarmass} and Theorem \ref{sphere}, we have $u$ is a standard bubble \;$\implies$\;$u=T_{g_{\mathbb{S}^n}}(u)$.
\item
The approach of this paper works also for Yamabe type problems on manifolds with boundary, Yamabe type problems of fractional and or higher order, their associated prescribed curvature problems, and their CR-version as well, see  \cite{alm1}, \cite{bc1}, \cite{bcch}, \cite{gam1}, \cite{gam2}, \cite{ja}, \cite{am},  \cite{gq1}, \cite{mw}, \cite{smw}, \cite{marquez1}, \cite{mayndai}, \cite{mayndia1}, \cite{nss} for more informations about those topics .
\item
The use of optimal control theory in curvature prescription problems from Conformal Geometry was initiated in our work\cite{ndiaye} regarding the $Q$-curvature equation. This work is a Yamabe version of \cite{ndiaye}.  
\item
For motivations to study geometric obstacle problems and optimal control problems in this geometric context, and the place of this work in the general control theory framework, see \cite{ndiaye}.
\item
Because of Proposition \ref{sameminimizer}, the optimal control problem \eqref{stateop} provides an other approach to the Yamabe minimizing problem.
\end{itemize}
\end{rem}
\vspace{7pt}

\noindent
To prove Theorem \ref{positivevarmass}-Theorem \ref{sphere}, we first use the variational characterization of the solution of conformal Laplacian obstacle problem \;$T_g(u)$ (see Lemma \ref{obstacleq})\;  to show that \;$T_g$\; is  idempotent, see Proposition \ref{nilpotent}. Next, using the idempotent property of \;$T_g$, we establish some monotonicity formulas fo $J^g$ and $I^g$ involving $T_g$, see Lemma \ref{decreasingformula}, Lemma \ref{minimaltunnel}, and Lemma \ref{decreasingformulaop}. Using the latter monotonicity formulas, we show that  any minimizer of \;$J^g$\; or any solution of the optimal control problem \eqref{stateop} is a fixed point of \;$T_g$, see Corollary \ref{rigidityminimizer} and Corollary \ref{rigidityminimizeri}. This allows us to show that  $\inf_{H^1_+(M, g)}I^g=\inf_{H^1_+(M, g)} J^g$ and \;$u\in H^1_+(M, g)$ is a minimizer of \;$I^g$\; is equivalent to $u$ is a mimimizer of $J^g$, see Proposition \ref{sameinf} and Proposition \ref{sameminimizer} . With this at hand, Theorem \ref{positivevarmass} follows from the works Yamabe\cite{yamabe}, Trudinger\cite{tru}, Aubin\cite{aubin} and Schoen\cite{schoen}, and the smooth regularity of positive solutions of the Yamabe equation \eqref{eq:sequation}. Furthermore, Theorem \ref{sphere} follows from the the fact that the standard bubbles are the only minimizers of $J^g$ on $H^1_+(M, g)$ when $(M, g)=(\mathbb{S}^n, g_{\mathbb{S}^n})$.
\vspace{7pt}

\noindent
The structure of the paper is as follows. In Section \ref{eq:notpre}, we collect some preliminaries and fix some notations. In Section \ref{opp}, we discuss the obstacle problem for the conformal Laplacian and show that \;$T_g$\; is  idempotent and homogeneous, and $I^g$ is scale invariant. In Section \ref{conformalrule}, we study the conformal transformation rules of  \;$T_g$\;and \;$I^g$. In Section \ref{monotonicity}, we establish some monotonicity formula for $J^g$ when passing from \;$u$\; to \;$T_g(u)$\; and present some applications. In Section \ref{monotonicityop}, we compare $I^g$ and $J^g$ , establish some monotonicity formula for $I^g$ when passing from \;$u$\; to \;$T_g(u)$, and present some applications. In Section \ref{comparingval}, we show $\inf_{H^1_+(M, g)}I^g=\inf_{H^1_+(M, g)} J^g$ and \;$u\in H^1_+(M, g)$ is a minimizer of \;$I^g$\; is equivalent to \;$u$\; is a minimizer of $J^g$. Moreover, we present the proof of Theorem \ref{positivevarmass}. Finally, in Section \ref{casesphere}, we derive a Sobolev type inequality involving the obstacle operator  \;$T_g$\; and present the proof of Theorem \ref{sphere}.
\section{Notations and Preliminaries}\label{eq:notpre}
 In this brief section, we fix our notations and give some preliminaries. First of all, from now until the end of the paper, \;$(M, g)$\; is the background $n$-dimensional closed Riemannian manifold with $n\geq 3$. Even if Theorem \ref{positivevarmass} and Theorem \ref{sphere} are stated with the background metric $g$, we will be working with some generality with $h\in [g]$. 
\vspace{7pt}
 
\noindent
For  \;$h\in [g]$, we recall the Yamabe functional \;$J^h$\; and its subcritical approximation \;$J_p^h$ ($1\leq p<2^*-1$)
 \begin{equation}\label{eq:defj1}
J^h(u):= \frac{\left<u,u\right>_h}{||u||^2_{L^{2^*}(M, h)}} ,\;\;\;\;u\in H^1_+(M, h),
\end{equation}
and
\begin{equation}\label{eq:defjt1}
J^h_p(u):= \frac{\left<u,u\right>_h}{||u||^2_{L^{p+1}(M, h)}} \;\;\;\;u\in H^1_+(M, h),
\end{equation}
where
\;$H^1_+(M, h)$, $\left<\cdot,\cdot\right>_h$, and  \;$L^{2^*}(M, h)$\; are as in Section \ref{intro} with \;$g$\; replaced \;$h$. Moreover, $L^{p+1}(M, h)$ is the standard Lebesgue space of functions which are $(p+1)$-integrable on \;$M$\; with respect to \;$g$ .
We set
\begin{equation}\label{valends}
J_{2^*-1}^h:=J^h.
\end{equation}
and have a family of functional \;$(J_p^h)_{1\leq p\leq 2^*-1}$ defined on \;$H^1_+(M, h)$. Clearly by definition of \;$J^h_p$, 
\begin{equation}\label{phjph}
J_p^h(\l u)=J_p^h(u),\;\; \l>0,\; u\in H^1_+(M, h).
\end{equation}
We recall that the conformal class $[g]$  of $g$ is
$$
[g]:=\{g_w:=w^{\frac{4}{n-2}}g, \;\;w\in C^{\infty}_+(M, g)\},
$$
and in this paper we use the notation $$g_{*}=*^{\frac{4}{n-2}}g.$$
\vspace{7pt}

\noindent
The following transformation rules are well-known and easy to verify (for \;$w\in C^{\infty}_+(M, g)$)
\begin{equation}\label{inspace}
wC^{\infty}_+(M, g_w)=C^{\infty}_+(M, g), 
\end{equation}
\begin{equation}\label{inspace}
wH^{1}_+(M, g_w)=H^1_+(M, g).
\end{equation}
\begin{equation}\label{invv}
dV_{g_w}=w^{2^*}dV_g
\end{equation}
\begin{equation}\label{invs}
\left<u, u\right>_{g_w}=\left<wu, wu\right>_{g}, \;\;\;u\in H^{1}_+(M, g_w)
\end{equation}
\begin{equation}\label{invdf}
||u||_{L^{2^*}(M, g_w)}=||uw||_{L^{2^*}(M, g)}, \;\;\;u\in H^{1}_+(M, g_w)
\end{equation}
and 
\begin{equation}\label{invyamabe}
J^{g_w}(u)=J^g(wu), \;\;\;u\in H^{1}_+(M, g_w).
\end{equation}
We recall also the fact  that \;$\inf_{H^1_+(M, g) }J^g$\; in the definition of \;$\mathcal{Y}(M, [g])$\; depends only on \;$[g]$\; can be seen from \eqref{inspace} and  \eqref{invyamabe}.
\vspace{7pt}

\noindent
For  \;$\mathcal{Y}(M, [g])>0$ and $h\in[g]$, we define the Yamabe optimal obstacle functional $I^h$ and its subcritical approximations  $I^h_p$ ($1\leq p<2^*-1$) by 
\begin{equation}\label{eq:defi1}
I^h(u):= \frac{\left<u,u\right>_h}{||T_h(u)||^2_{L^{2^*}(M, h)}} ,\;\;\;\;u\in H^1_+(M, h),
\end{equation}
and
\begin{equation}\label{eq:defit1}
I^h_p(u):= \frac{\left<u,u\right>_h}{||T_h(u)||^2_{L^{p+1}(M, h)}} \;\;\;\;u\in H^1_+(M, h),
\end{equation}
where  \;$T_h$\; is as in \eqref{defoo} with \;$g$\; replaced by \;$h$.
We set
\begin{equation}\label{defend}
I_{2^*-1}^h:=I^h
\end{equation}
and have a family of functional \;$(I_p^h)_{1\leq p\leq 2^*-1}$\; defined on \;$H^1_+(M, h)$. 
\vspace{7pt}

\noindent
When $\mathcal{Y}(M, [g])>0$, the following Sobolev type inequality holds.
\begin{lem}\label{sobolev1}
Assuming \;$\mathcal{Y}(M, [g])>0$ and $h\in[g]$, then for $u\in H^1_+(M, h)$
$$
||u||_{L^{2*}(M, h)}\leq \frac{1}{\sqrt{\mathcal{Y}(M, [g])}}||u||_h
$$
\end{lem}
\vspace{7pt}

\noindent
When  \;$(M, g)=(\mathbb{S}^n, g_{\mathbb{S}^n})$, we have the following well-known stronger version of the latter Sobolev inequality.
\begin{lem}\label{sobolev2}
Assuming \;$(M, g)=(\mathbb{S}^n, g_{\mathbb{S}^n})$ and $h=g_w\in[g]$, then for  $u\in H^1_+(M, h)$,
\begin{equation}\label{inequality0}
||u||_{L^{2*}(M, h)}\leq \frac{1}{\sqrt{\mathcal{Y}(\mathbb{S}^n, [g_{\mathbb{S}^n}])}}||u||_h,
\end{equation}
with equality holds if and only if $$u\in C^{\infty}_+(M, h)\;\;\text{and}\;\;R_{g_{wu}}=const.$$
\end{lem}
\begin{rem}
We recall that the explicit value of  \;$\mathcal{Y}(\mathbb{S}^n, [g_{\mathbb{S}^n}])$\; is known.
\end{rem}

\section{Obstacle problem for the conformal Laplacian}\label{opp}
In this section, we study the obstacle problem for the conformal Laplacian $L_h$ with \;$h\in [g]$\ under the assumption $\mathcal{Y}(M, [g])>0$. Indeed in analogy to the classical obstacle problem for the Laplacian $\D_h$, given \;$u\in H^1_+(M, h)$, we look for a solution to the minimization problem 
\begin{equation}\label{obspan}
\min_{v\in H^1_+(M, h),\\\;\;v\geq u}\left<v, v \right>_h.
\end{equation} 
We anticipate that the assumption $\mathcal{Y}(M, [g])>0$ guarantees that the minimization problem \eqref{obspan} is well posed and $$||v||_h:=\sqrt{\left<v, v\right>_{h}}$$ is a norm, since \; $L_h\geq 0$\; and \;$\ker L_h=\{0\}$.
\vspace{7pt}

\noindent
We start with the following lemma providing the existence and unicity of solution for the obstacle problem for the conformal Laplacian \;$L_h$\; \eqref{obspan}. 
\begin{lem}\label{obstacleq}
Assuming that \;$\mathcal{Y}(M, [g])>0$ and $h\in[g]$, then  for \;$u\in H^1_+(M, h)$,
there exists a unique $T_h(u)\in H^1_+(M, h)$ such that
\begin{equation}\label{tg}
||T_h(u)||^2_h=\min_{v\in H^1_+(M, h),\\\;\;v\geq u}||v||^2_h
\end{equation}
\end{lem}
\begin{pf}
Since \;$\mathcal{Y}(M, [g])>0$, then  $L_h\geq 0$\; and \;$\ker L_h=\{0\}$. Thus \;$<\cdot, \cdot>_h$\; defines a scalar product on \;$H^1_+(M, h)$\; inducing a norm $||\cdot||_h$ equivalent to the standard \;$H^1(M, h)$-norm on \;$H^1_+(M, h)$. Hence, as in the classical obstacle problem for the Laplacian $\D_h$, the lemma follows from Direct Methods in the Calculus of Variations.
\end{pf}
\vspace{7pt}

\noindent
We study now some properties of the obstacle solution map (or state map) \;$T_h\; : H^1_+(M, h)\longrightarrow  H^1_+(M, h)$. We start with the following algebraic one.
\begin{pro}\label{nilpotent}
Assuming that  \;$\mathcal{Y}(M, [g])>0$\; and  \;$h\in[g]$, then the state map \;$T_h\; : H^1_+(M, h)\longrightarrow  H^1_+(M, h)$\; is idempotent, i.e $$T^2_h=T_h.$$
\end{pro}

\begin{pf}
Let $v\in H^1_+(M, h)$ such that $v\geq T_h(u)$. Then $T_h(u)\geq u$ implies
$v\geq u$. Thus by minimality, \;$$||v||_h \geq ||T_h(u)||_h.$$ Hence, since  \;$T_h(u)\in H^1_+(M, h)$\; and \;$T_h(u)\geq T_h(u)$\; then by uniqueness 
$$
T_h(T_h(u))=T_h(u),
$$
thereby ending the proof.
\end{pf}
\vspace{7pt}

\noindent
We have the following lemma showing that $T_h$ is positively homogeneous.
\begin{lem}\label{pht}
Assuming  \;$\mathcal{Y}(M, [g])>0$\; and  \;$h\in [g]$, then for  \;$\l>0$, 
$$
T_h(\lambda u)=\lambda T_h(u), \;\;\;\forall u\in H^1_+(M, h).
$$
\end{lem}
\begin{pf}
Let $\l u\leq v\in H^1_+(M, h)$. Then $\l>0$ implies $\l^{-1}v\geq u$. Thus, since \;$\l^{-1}v\in H^1_+(M, h)$, then by minimality 
$$
||\l^{-1}v||_h\geq ||T_h(u)||_h.
$$
So by positive homogeneity of \; $||\cdot||_h$, we get
$$
||v||_h\geq ||\l T_h(u)||_h
$$
Hence, since $\l u\leq \l T_h(u)\in H^1_+(M, h)$, then by uniqueness we obtain
$$
T_h(\l u)=\l T_h(u),
$$
as desired.
\end{pf}
\vspace{7pt}

\noindent
Lemma \ref{pht} implies the following analogue of formula \eqref{phjph} for $I^h_p$.
\begin{cor}\label{scaleinv}
Assuming  \;$\mathcal{Y}(M, [g])>0$,  \;$h\in [g]$, and \;$1\leq p\leq 2^*-1$, then for \;$\l>0$,
$$
I^h_p(\lambda u)=I_p^h(u), \;\;\;\forall u\in H^1_+(M, h).
$$
\end{cor}
\begin{pf}
It follows directly from the definition of \;$I^h_p$ (see \eqref{eq:defi1}-\eqref{defend}), Lemma \ref{pht}), and the positive homogeneity of norms. Indeed, \eqref{eq:defi1}-\eqref{defend} imply
\begin{equation}\label{invs1}
I^h_p(\lambda u)=\frac{||\l u||^2_h}{||T_h(\l u)||^2_{L^{p+1}(M, h)}}. 
\end{equation}
Thus, Lemma \ref{pht} and \eqref{invs1} give
\begin{equation}\label{invs2}
I^h_p(\lambda u)=\frac{||\l u||^2_h}{||\l T_h(u)||^2_{L^{p+1}(M, h)}}.
\end{equation}
So, the positive homogeneity of norms and \eqref{invs2} imply
\begin{equation}\label{invs3}
I^h_p(\lambda u)=\frac{||u||^2_h}{||T_h(u)||^2_{L^{p+1}(M, h)}} .
\end{equation}
Hence, using again \eqref{eq:defi1}-\eqref{defend} combined with \eqref{invs3}, we get
$$
I^h_p(\lambda u)=I^h_p(u),
$$
as desired.
\end{pf}
\begin{rem}
We would like to point out that Lemma \ref{pht}  and Corollary \ref{scaleinv} are not used in the proof of Theorem \ref{positivevarmass} and Theorem \ref{sphere}. We decide to put them in this paper, because beside of being of independent interest, Lemma \ref{scaleinv} is a direct analogue \eqref{phjph} and is useful for a direct variational study of \;$I^h_p$\; using the Direct Methods in the Calculus of Variations. Moreover, Lemma \ref{pht} is needed in the proof of Corollary \ref{scaleinv} .
\end{rem}
\section{Transformation rules of  \;$T_{h}$\; and  \;$I^{h}$\; for \;$h\in [g]$}\label{conformalrule}
In this section, we study the transformation rules of \;$T_h$\; when \;$h$\; varies in \;$[g]$\; and use it to establish an analogue of formula \eqref{invyamabe} for the Yamabe optimal control function $I^h$. 
\vspace{7pt}

\noindent
We adopt the notation \;$h=g_w=w^{\frac{4}{n-2}}g$\; with \;$w\in C^{\infty}_+(M, g)$\; and have:
\begin{lem}\label{tth}
Assuming that  \;$\mathcal{Y}(M, [g])>0$\; and \;$w\in C^{\infty}_+(M, g)$, then
$$
T_{g_w}(u)=w^{-1}T_g(wu), \;\;\;u\in H^1_+(M, g_w).
$$
\end{lem}
\begin{pf}
Let \;$v\in H^1_+(M, g_w)$\; with \;$v\geq u$. Then using formula \eqref{invs}, we get
$$
||v||_{g_w}=||wv||_g.
$$
Thus, since \;$v\geq u$\; and \;$w>0$\; implies \;$wv\geq uw$, then  using $wu, \;wv\in wH^1_+(M, g_w)$ and  \eqref{inspace}, we have by minimality
$$
||v||_{g_w}\geq ||T_g(wu)||_{g}.
$$
So, using again formula \eqref{invs}, we obtain
$$
||v||_{g_w}\geq ||w^{-1}T_g(wu)||_{g_w}
$$
Now, since  $w^{-1}T_g(wu)\geq w^{-1}wu=u$ and $w^{-1}T_g(wu)\in H^1_+(M, g_w)$ (see \eqref{inspace}), then by uniqueness 
$$
T_{g_w}(u)=w^{-1}T_g(wu),
$$ 
thereby ending the proof.
\end{pf}
\vspace{7pt}

\noindent
As a consequence of Lemma \ref{tth}, we have:
\begin{cor}\label{invfixpt}
Assuming that  \;$\mathcal{Y}(M, [g])>0$\; and  \;$w\in C^{\infty}_+(M, g)$, then
$$
Fix(T_{g_w})=w^{-1}Fix(T_g),
$$
where for $h\in [g]$,  
\begin{equation}\label{deffix}
Fix(T_h):=\{u\in H^1(M, h):\;\;T_h(u)=u\}.
\end{equation}
\end{cor}
\begin{pf}
 Lemma \ref{tth} and \eqref{deffix} imply
 $$
u\in Fix(T_{g_w})\Longleftrightarrow\;\;T_{g_w}(u)=u\;\;\Longleftrightarrow\;\; w^{-1}T_g(wu)=u.
$$ 
Thus
$$u\in Fix(T_{g_w})\;\;\Longleftrightarrow\;\;T_g(wu)=wu\;\;\Longleftrightarrow \;\;wu\in Fix(T_g).$$
Hence$$
u\in Fix(T_{g_w})\;\;\;\Longleftrightarrow\;\;\;u\in w^{-1}Fix(T_g),
$$
as desired
\end{pf}
\vspace{7pt}

\noindent
Lemma \ref{tth} also implies the following corollary.
\begin{cor}\label{invd}
Assuming that  \;$\mathcal{Y}(M, [g])>0$\; and \;$w\in C^{\infty}_+(M, g)$, then
$$
||T_{g_w}(u)||_{L^{2^*}(M, g_w)}=||T_{g}(wu)||_{L^{2^*}(M, g)}, \;\;\;u\in H^1_+(M, g_w).
$$
\end{cor}
\begin{pf}
By definition of the \;$L^{2^*}(M, g_w)$-norm, we have
$$
||T_{g_w}(u)||_{L^{2^*}(M, g_w)}^{2^*}=\int_MT_{g_w}(u)^{2^*}dV_{g_w}.
$$
Thus, using \eqref{invv} and Lemma \ref{tth}, we get
$$
||T_{g_w}(u)||_{L^{2^*}(M, g_w)}^{2^*}=\int_MT_{g}(wu)^{2^*}dV_{g}
$$
So, using the definition of the \;$L^{2^*}(M, g)$-norm, we obtain
$$
||T_{g_w}(u)||_{L^{2^*}(M, g_w)}^{2^*}=||T_{g}(wu)||_{L^{2^*}(M, g)}^{2^*}.
$$
Hence, we get 
$$
||T_{g_w}(u)||_{L^{2^*}(M, g_w)}=||T_{g}(wu)||_{L^{2^*}(M, g)},
$$
as desired.
\end{pf}
\vspace{7pt}

\noindent
As a consequence of Corollary \ref{invd}, we have the following analogue of formula \eqref{invyamabe} for the Yamabe optimal obstacle functional.
\begin{cor}\label{invi}
Assuming that  \;$\mathcal{Y}(M, [g])>0$\; and \;$w\in C^{\infty}_+(M, g)$, then for  \;\;$u\in H^1_+(M, g_w)$, 
$$
I^{g_w}(u)=I^g(wu).
$$
\end{cor}
\begin{pf}
By definition of \;$I^{g_w}$ (see \eqref{eq:defi1}), we have
$$
I^{g_w}(u)=\frac{||u||^2_{g_w}}{||T_{g_w}(u)||^2_{L^{2^*}(M, g_w)}}
$$
Thus, using \eqref{invs} and Corollary \ref{invd}, we obtain
$$
I^{g_w}(u)=\frac{||wu||^2_{g}}{||T_g(wu)||^2_{L^{2^*}(M, g)}}
$$
Hence, using the definition of \;$I^g$ (see \eqref{eq:defi1}), we get
$$
I^{g_w}(u)=I^g(wu),
$$
as desired.	
\end{pf}
\vspace{5pt}

\noindent
Corollary \ref{invi} implies.
\begin{cor}\label{cor14}
Assuming that  \;$\mathcal{Y}(M, [g])>0$\; and \;$w\in C^{\infty}_+(M, g)$, then
$$
\inf_ {u\in H^1_+(M, g_w)} I^{g_w}(u)=\inf _ {u\in H^1_+(M, g)} I^g(u).
$$
\end{cor}
\begin{pf}
 \eqref{inspace} implies
 $$
 \inf _ {u\in H^1_+(M, g)} I^g(u)=\inf _ {u\in wH^1_+(M, g_w)} I^g(u)
 $$
 Thus, setting $u=w\bar u$, we get
 $$
 \inf _ {u\in H^1_+(M, g)} I^g(u)=\inf _ {\bar u\in H^1_+(M, g_w)} I^g(w\bar u)
 $$
 Hence, Corollary \ref{invi} implies
 $$
 \inf _ {u\in H^1_+(M, g)} I^g(u)=\inf _ {\bar u\in H^1_+(M, g_w)} I^{g_w}(\bar u),
 $$
as desired.
 \end{pf}
\vspace{7pt}

\noindent
Similar to the Yamabe invariant \;$\mathcal{Y}(M, [g])$, Corollary \ref{cor14} justifies the following definition.
\begin{df}\label{defyp}
Assuming that  \;$\mathcal{Y}(M, [g])>0$ and \;$h\in [g]$, then
$$
\mathcal{Y}_{oc}(M, [h]):=\inf _ {u\in H^1_+(M, h)} I^h(u).
$$
\end{df}
\vspace{7pt}

\noindent
\begin{rem}
Clearly \;for $h\in [g]$, $$\mathcal{Y}_{oc}(M, [h])=\mathcal{Y}_{oc}(M, [g]).$$
\end{rem}


\vspace{7pt}

\noindent
In the work of Yamabe\cite{yamabe}, the following family of real numbers (not conformal invariant) was introduced
$$
\mathcal{Y}^p(M, h):=\inf _ {u\in H^1_+(M, h)} J^h_p(u), \;\;h\in[g], \;1\leq p<2^*-1.
$$
We set
$$
\mathcal{Y}^{2^*-1}(M, h)=\mathcal{Y}(M, [h])
$$
and have a family of real numbers $(\mathcal{Y}^p(M, h))_{1\leq p\leq 2^*-1}$ which is conformally invariant for $p=2^*-1$.
\vspace{7pt}

\noindent
Similarly, for $\mathcal{Y}(M, [g])>0$ we define
$$
\mathcal{Y}^p_{oc}(M, h):=\inf _ {u\in H^1_+(M, h)} I^h_p(u), \;\;h\in[g], \;1\leq p<2^*-1
$$
and
$$
\mathcal{Y}^{2^*-1}_{oc}(M, h)=\mathcal{Y}_{oc}(M, [h]).
$$
This defines a family of real numbers $(\mathcal{Y}^p_{oc}(M, h))_{1\leq p\leq 2^*-1}$ which is conformally invariant for $p=2^*-1$
\section{Monotonicity formula for \;$J^h_p$, $h\in [g]$, and $1\leq p\leq 2^*-1$}\label{monotonicity}
In this section, we present a monotonicity formula for \;$J^h_p$\; when passing from \;$u$\; to \;$T_h(u)$\; for \;$h\in[g]$\; with \;$\mathcal{Y}(M, [g])>0$ \; and \;$1\leq p\leq 2^*-1$. Moreover, we present some applications on the relation between the ground state of \;$J^h_p$\; and the the fixed points of \;$T_h$.  
\vspace{7pt}

\noindent
The monotonicity formula reads as follows.
\begin{lem}\label{decreasingformula}
Assuming that  \;$\mathcal{Y}(M, [g])>0$, \;$h\in[g]$\; and \;$1\leq p\leq 2^{*}-1$, then for \;$u\in H^1_+(M, h)$, $$
J_p^h(u)-J_p^h(T_h(u))\geq \frac{1}{||T_h(u)||^2_{L^{p+1}(M, h)}}\left[||u||^2_h-||T_h(u))||^2_h\right]\geq 0.
$$
\end{lem}
\begin{pf}
Using the definition of \;$J_p^h$ (see \eqref{eq:defj1}-\eqref{valends}), we have
\begin{equation}
J_p^h(u)-J_p^h(T_h(u))= \frac{||u||^2_h}{||u||^2_{L^{p+1}(M, h)}}-\frac{||T_h(u)||^2_h}{||T_h^2(u)||^2_{L^{p+1}(M, h)}}.
\end{equation}
Thus, using \;$T_h^2(u)=T_h(u)$\; (see Proposition \ref{nilpotent}), we get
\begin{equation}
J_p^h(u)-J_p^h(T_h(u))= \frac{||u||^2_h}{||u||^2_{L^{p+1}(M, h)}}-\frac{||T_h(u)||^2_h}{||T_h(u)||^2_{L^{p+1}(M, h)}}.
\end{equation}
Hence the result follows from  \;$T_h(u)\geq u> 0$, and and the definition of \;$T_h$\; (see Lemma \ref{obstacleq}).
\end{pf}
\vspace{7pt}

\noindent
Lemma \ref{decreasingformula} implies the following rigidity result.
\begin{cor}\label{rigidity}
Assuming that  \;$\mathcal{Y}(M, [g])>0$, \;$h\in[g]$, and \;$1\leq p\leq 2^{*}-1$, then for\;$u\in H^1_+(M, h)$,
\begin{equation}\label{ineq}
J_p^h(T_h(u)) \leq  J_p^h(u)
\end{equation}
and  
\begin{equation}\label{eql}
J_p^h(u)=J_p^h(T_h(u))\implies u \in Fix(T_h).
\end{equation}
\end{cor}
\begin{pf}
Using Lemma \ref{decreasingformula}, we have
\begin{equation}\label{eq1}
J_p^h(u)-J_p^h(T_h(u))\geq  \frac{1}{||T_h(u)||^2_{p+1}}\left[||u||^2_h-||T_h(u)||^2_h\right]\geq 0.
\end{equation}
Thus, \eqref{ineq} follows from \eqref{eq1}. If \;$J_p^h(u)=J_p^h(T_h(u))$, then \eqref{eq1} implies
$$
||u||^2_h=||T_h(u)||^2_h.
$$
Hence, since \;$u\geq u\in H^1_+(M, h)$, then the uniqueness part in Lemma \ref{obstacleq} implies
$$
u=T_h(u).
$$
Hence, using \eqref{deffix}, we have
$$
u\in Fix(T_h),
$$
thereby ending the proof of the corollary.
\end{pf}


\vspace{7pt}

\noindent
Corollary \ref{rigidity} implies that minimizers of  \;$J_p^h$\; on \;$H^1_+(M, h)$\; belongs to \;$Fix(T_h)$. Indeed, we have:
\begin{cor}\label{rigidityminimizer}
Assuming that  \;$\mathcal{Y}(M, [g])>0$,  $h\in[g]$, and  \;$1\leq p\leq 2^{*}-1$, then for \;$u\in H^1_+(M, h)$, $$J_p^h(u)=\mathcal{Y}^p(M , h)\implies u\in Fix(T_h).$$
\end{cor}
\begin{pf}
$J_p^h(u)=\mathcal{Y}^{p}(M , h)$\; implies
\begin{equation}\label{eq3}
J_p^h(u)\leq J_p^h(T_h(u)).
\end{equation}
Thus combining \eqref{ineq} and \eqref{eq3}, we get
\begin{equation}\label{eq4}
J_p^h(u)= J_p^h(T_h(u)).
\end{equation}
Hence, combining \eqref{eql} and \eqref{eq4}, we obtain
$$
u\in Fix (T_h),
$$
as desired
\end{pf}

\vspace{7pt}

\noindent
\begin{rem}\label{rem2022}
Under the assumption of Corollary \ref{rigidity}, we have Proposition \ref{nilpotent} and Corollary \ref{rigidity} imply that we can assume without loss of generality that any minimizing sequence \;$(u_l)_{l\geq 1}$\; of \;$J_p^h$\; on  \;$H^1_+(M, h)$\;satisfies
$$
u_l\in Fix(T_h), \;\;\;\;\forall l\ge 1.
$$
Indeed, suppose \;$u_l$\; is a minimizing sequence for \;$J_p^h$\; on \;$H^1_+(M, h)$. Then \;$u_l\in H^1_+(M, h)$\; and 
$$
J_p^h(u_l)\longrightarrow \inf_{H^1_+(M, h)}J_p^h.
$$
Thus by definition of infimum and Corollary \ref{rigidity}, we  have
$$
\inf_{H^1_+(M, h)}J_p^h\leq J_p^h(T_h(u_l))\leq J_p^h(u_l).
$$
This implies $$J_p^h(T_h(u_l))\longrightarrow \inf_{H^1_+(M, h)}J_p^h.$$ Hence setting $$\hat u_l=T_h(u_l),$$ and using Proposition \ref{nilpotent}, we get
$$J_p^h(\hat u_l)\longrightarrow \inf_{H^1_+(M, h)}J_p^h\;\;\; \;\;\text{and}\;\; \;\;\;\;\hat u_l=T_g(\hat u_l)$$ as desired.
\end{rem}

\section{Monotonicity formula for $I^h_p$, $h\in [g]$,  $1\leq p\leq 2^*-1$}\label{monotonicityop}
In this section, we derive a monotonicity formula for $I_p^h$ similar to the one for $J_p^h$ derived in the previous section for $h\in [g]$\; with \;$\mathcal{Y}(M, [g])>0$\; and \;$1\leq p\leq 2^{*}-1$. Moreover, we present some applications for \;$I^h_p$\; similar to the ones done for $J^h_p$ in the previous section.
\vspace{7pt}

\noindent
We have the following monotonicity formula for the \;$p$-Yamabe optimal control functional \;$I_p^h$.
\begin{lem}\label{decreasingformulaop}
Assuming that  \;$\mathcal{Y}(M, [g])>0$, $h\in [g]$, and  \;$1\leq p\leq 2^{*}-1$, then for  \;$u\in H^1_+(M, h)$,
$$
I_p^h(u)-I_p^h(T_h(u))=\frac{1}{||T_h(u)||^2_{L^{p+1}(M, h)}}\left[||u||^2_h-||T_h(u)||^2_h\right]\geq 0.
$$
\end{lem}
\begin{pf}
By definition of \;$I_p^h$ (see \eqref{eq:defi1}-\eqref{defend}), we have
$$
I_p^h(u)-I_p^h(T_h(u))=\frac{||u||^2_h}{||T_h(u)||^2_{L^{p+1}(M, h)}}-\frac{||T_g(u)||^2_h}{||T_h^2(u)||^2_{L^{p+1}(M, h)}}.
$$
Using \;$T_h^2(u)=T_h(u)$\;(see Proposition \ref{nilpotent}) and the definition of \;$T_h$\; (see Lemma \ref{obstacleq}), we  get
$$I_p^h(u)-I_p^h(T_h(u))=\frac{1}{||T_h(u)||^2_{L^{p+1}(M, h)}}\left[||u||^2_h-||T_g(u)||^2_h\right]\geq 0,$$
thereby ending the proof.
\end{pf}
\vspace{7pt}

\noindent
Similar to the previous section,  Lemma \ref{decreasingformulaop} implies the following rigidity result.
\begin{cor}\label{rigidityop}
Assuming that \;$\mathcal{Y}(M, [g])>0$,  $h\in[g]$, and  \;$1\leq p\leq 2^{*}-1$, then for  \;$u\in H^1_+(M, h)$,
\begin{equation}\label{ineqop}
I_p^h(T_h(u)) \leq  I_p^h(u)
\end{equation}
and  
\begin{equation}\label{eqlop}
I_p^h(u)=I_p^h(T_h(u))\implies u \in Fix(T_h).
\end{equation}
\end{cor}
\begin{pf}
Using Lemma \ref{decreasingformulaop}, we have
\begin{equation}\label{eq1op}
I_p^h(u)-I_p^h(T_h(u))=\frac{1}{||T_h(u)||^2_{p+1}}\left[||u||^2_h-||T_h(u)||^2_h\right]\geq 0.
\end{equation}
Thus, \eqref{ineqop} follows from \eqref{eq1op}. If \;$I_p^h(u)=I_p^h(T_h(u))$, then \eqref{eq1op} implies
$$
||u||^2_h=||T_h(u)||^2_h.
$$
Hence, since \;$u\geq u\in H^1_+(M, h)$, then  as in the previous section the uniqueness part in Lemma \ref{obstacleq} implies
$$
u\in Fix(T_h),
$$
thereby ending the proof of the corollary.
\end{pf}
\vspace{7pt}

\noindent
As in the previous section, corollary \ref{rigidity} imply that minimizers of  \;$I_p^h$\; belongs to  \;$Fix(T_h)$.
\begin{cor}\label{rigidityminimizeri}
Assuming that  \;$\mathcal{Y}(M, [g])>0$,  $h\in[g]$\;  and  \;$1\leq p\leq 2^{*}-1$, then for \;$u\in H^1_+(M, h)$, $$I_p^h(u)=\mathcal{Y}_{oc}^p(M, h)\implies u \in Fix(T_h).$$
\end{cor}
\begin{pf}
$I_p^h(u)=\mathcal{Y}^p_{oc}(M, h)$\; implies
\begin{equation}\label{eq3op}
I_p^h(u)\leq I_p^h(T_h(u)).
\end{equation}
Thus combining \eqref{ineqop} and \eqref{eq3op}, we get
\begin{equation}\label{eq4op}
I_p^h(u)= I_p^h(T_h(u)).
\end{equation}
Hence, combining \eqref{eqlop} and \eqref{eq4op}, we obtain
$$
u\in Fix (T_h),
$$
thereby ending the proof.

\end{pf}
\vspace{7pt}

\noindent
\begin{rem}
Under the assumptions of Corollary \ref{decreasingformulaop} and using the same argument as in Remark \ref{rem2022}, we have that Proposition \ref{nilpotent} and Corollary \ref{decreasingformulaop} imply that for a minimizing sequence \;$(u_l)_{l\geq 1} $\; of \;$I_p^h$\; on  \;$H^1_+(M, h)$, we can  assume without loss of generality that
$$
u_l\in Fix(T_h), \;\;\;\;\forall l\ge 1.
$$
\end{rem}
\vspace{10pt}

\noindent
\section{Comparing \;$\mathcal{Y}^p(M, h)$\; and \;$\mathcal{Y}_{oc}^p(M, h)$, $h\in[g]$, and \;$1\leq p\leq 2^*-1$}\label{comparingval}
In this section, for \;$\mathcal{Y}(M, [g])>0$,  \;$h\in [g]$,  \;$1\leq p\leq 2^*-1$, we show  \;$\mathcal{Y}^p(M, h)=\mathcal{Y}_{oc}^p(M, h)$\; and \;$J^h_p(u)=\mathcal{Y}^p(M, h)\Longleftrightarrow I^h_p(u)=\mathcal{Y}_{oc}^p(M, h)$. As a consequence, we deduce Theorem \ref{positivevarmass}. 
\vspace{7pt}

\noindent
We start with the following comparison result showing that  \;$I_p^h\leq J_p^h$\; and \;$I_p^h=I_p^h$\; on \;$T_h(H^1_+(M, h))$\; the range of \;$T_h$. 
\begin{lem}\label{minimaltunnel}
Assuming that \;$\mathcal{Y}(M, [g])>0$,  $h\in[g]$\;  and  \;$1\leq p\leq 2^{*}-1$, then 
\begin{equation}\label{lesval}
I_p^h\leq J_p^h\;\;\;\;\;\;\text{on}\;\;\;\;\; H^1_+(M, h)
\end{equation}
and
\begin{equation}\label{sameval}
J_p^h\circ T_h=I_p^h\circ T_h\;\;\;\;\text{on}\;\;\; H^1_+(M, h).
\end{equation}
\end{lem}
\begin{pf}
By definition of \;$J_p^h$\; and \;$I_p^h$ (see \eqref{eq:defj1}-\eqref{valends}, and \eqref{eq:defi1})-\eqref{defend}, for $u\in H^1_+(M, h)$ we have
\begin{equation}\label{app1}
J_p^h(u)-I_p^h(u)=\frac{||u||^2_h}{||u||^2_{L^{p+1}(M, h)}}-\frac{||u||^2_h}{||T_h(u)||^2_{L^{p+1}(M, h)}}.
\end{equation}
Thus \eqref{lesval} follows from \;$T_h(u)\geq u>0$ and \eqref{app1}. Moreover, we have
$$
J_p^h(T_h(u))-I_p^h(T_h(u))= \frac{||T_h(u)||^2_h}{||T_h(u)||^2_{L^{p+1}(M, h)}}-\frac{||T_h(u)||^2_h}{||T_h^2(u)||^2_{L^{p+1}(M, h)}}.
$$
Hence, \;$T_h^2(u)=T_h(u)$\; (see Proposition \ref{nilpotent}) implies \;$$J_p(T_h(u))=I_p^h(T_h(u)),$$
thereby ending the proof.
\end{pf}
\begin{rem}
Clearly \;$T_h(u)\geq u>0$ and \eqref{app1} imply
$$
I_p^h(u)=J_p^h(u)\Longleftrightarrow u=T_h(u)
$$
\end{rem}
\vspace{7pt}

\noindent
Corollary \ref{rigidity}, Corollary \ref{rigidityop}, and Lemma \ref{minimaltunnel} imply \;$\mathcal{Y}^p(M, h)=\mathcal{Y}^p_{oc}(M, h)$. Indeed, we have:
\begin{pro}\label{sameinf}
Assuming that \;$\mathcal{Y}(M, [g])>0$,  \;$h\in[g]$\;  and  \;$1\leq p\leq 2^{*}-1$, then
$$
\mathcal{Y}^p(M, h)=\mathcal{Y}^p_{oc}(M, h).
$$
\end{pro}
\begin{pf}
Using \eqref{ineq}, we get
\begin{equation}\label{relinf1}
\mathcal{Y}^p(M, h)=\inf_{H^1_+(M, h)}J_p^h\circ T_h.
\end{equation}
Similarly, \eqref{ineqop} implies
\begin{equation}\label{relinf2}
\mathcal{Y}^p_{oc}(M, h)=\inf_{H^1_+(M, h)}I_p^h\circ T_h.
\end{equation}
Hence, the result follows from \eqref{sameval}, \eqref{relinf1}, and \eqref{relinf2}.
\end{pf}
\vspace{7pt}

\noindent
Combining Corollary \ref{rigidityminimizer}, Corollary \ref{rigidityminimizeri}, \eqref{sameval} and Proposition \ref{sameinf}, we have.
\begin{pro}\label{sameminimizer}
Assuming that \;$\mathcal{Y}(M, [g])>0$,  $h\in[g]$\;  and  \;$1\leq p\leq 2^{*}-1$, then for \;$u\in H^1_+(M, h)$, 
$$
J^h_p(u)=\mathcal{Y}^p(M, h)\;\;\;\Longleftrightarrow\;\; I^h_p(u)=\mathcal{Y}^p_{oc}(M, h).
$$
\end{pro}
\begin{pf}
If $J^h_p(u)=\mathcal{Y}^p(M, h)$, then Corollary \ref{rigidityminimizer} implies
$$
u\in Fix(T_h).
$$
So, using \eqref{deffix}, we get
$$
u=T_h(u)
$$
Thus, \eqref{sameval} gives
$$
I^h_p(u)=J^h_p(u).
$$
Hence, using Proposition \ref{sameinf}, we get $$I^h_p(u)=\mathcal{Y}^p_{oc}(M, h).$$
If $I^h_p(u)=\mathcal{Y}^p_{oc}(M, h)$, then  similarly as before, Corollary \ref{rigidityminimizeri} implies
$$
u=T_h(u).
$$
Thus, \eqref{sameval} gives
$$
J^h_p(u)=I^h_p(u).
$$
Hence, using Proposition \ref{sameinf}, we get $$J^h_p(u)=\mathcal{Y}^p(M, h),$$ and this completes the proof.
\end{pf}
\vspace{7pt}

\noindent
\begin{pfn}{ of Theorem \ref{positivevarmass}}\\
It follows from the works of Yamabe\cite{yamabe}, Trudinger \cite{tru}, Aubin \cite{aubin1}, and Schoen\cite{schoen}, Corollary \ref{rigidityminimizer}, Proposition \ref{sameminimizer}, and the smooth regularity of positive solutions of  the Yamabe equation \eqref{eq:sequation}. Indeed, for point 1), Proposition \ref{sameminimizer} implies
\begin{equation}\label{imp00}
J^g(u)=\mathcal{Y}(M, [g]).
\end{equation}
Thus as a critical point of $J^g$ and the  smooth regularity of positive solutions of the Yamabe equation \eqref{eq:sequation}
\begin{equation}\label{imp01}
u\in C^{\infty}_+(M, g)\;\;\;\text{ and}\;\;\;R_{g_u}=const>0.
\end{equation}
Corollary  \ref{rigidityminimizer}, \eqref{deffix}, and \eqref{imp00} imply
\begin{equation}\label{imp02}
u\in Fix(T_g)\Longleftrightarrow u=T_g(u).
\end{equation}
Hence, point 1) follows from \eqref{imp01} and \eqref{imp02}.\\
For point 2), by the works of Yamabe\cite{yamabe}, Trudinger \cite{tru}, Aubin \cite{aubin1}, and Schoen\cite{schoen}, Corollary \ref{rigidityminimizeri}, and the smooth regularity of positive solutions of  the Yamabe equation \eqref{eq:sequation}, there exists $u_{\min}\in C^{\infty}_+(M, g)$ such that 
\begin{equation}\label{imp1}
J^g(u_{min})=\mathcal{Y}(M, [g])
\end{equation}
and
\begin{equation}\label{imp2}
R_{g_{u_{min}}}=\mathcal{Y}(M, [g]).
\end{equation}
Thus Proposition \ref{sameminimizer} and \eqref{imp1} imply
\begin{equation}\label{imp3}
I^g(u_{min})=\mathcal{Y}_{oc}(M, [g]).
\end{equation}
Hence, point 2) follows from \eqref{imp2} and \eqref{imp3}.
\end{pfn}

\section{Obstacle problem and Sobolev type inequality}\label{casesphere}
In this section, we discuss some Sobolev type inequalities related to the obstacle problem for the conformal Laplacian. In particular, we specialize to the case of the \;$n$-dimensional standard sphere \;$(\mathbb{S}^n, g_{\mathbb{S}^n})$.
\vspace{7pt}

\noindent
Lemma \ref{sobolev1} (or Proposition \ref{sameinf}) implies the following obstacle Sobolev type inequality.
\begin{lem}\label{sobolev1op}
Assuming $\mathcal{Y}(M, [g])>0$\; and \;$h\in[g]$, then for \;$u\in H^1_+(M, h)$
$$
|||T_h(u)||_{L^{2*}(M, h)}\leq \frac{1}{\sqrt{\mathcal{Y}(M, [g])}}||u||_h
$$
\end{lem}
\begin{pf}
Lemma \ref{sobolev1} implies
$$
|||T_h(u)||_{L^{2*}(M, h)}\leq \frac{1}{\sqrt{\mathcal{Y}(M, [g])}}||T_h(u)||_h
$$
From the minimality property of $T_h(u)$, we get
$$
||T_h(u)||_h\leq ||u||_h.
$$
Thus combining the two inequalities, we obtain
$$
|||T_h(u)||_{L^{2*}(M, h)}\leq \frac{1}{\sqrt{\mathcal{Y}(M, [g])}}||u||_h
$$
as desired. On the other hand this can be deduced also using Proposition \ref{sameinf}. Indeed, using the definition of \;$I^h$\; and \;$\mathcal{Y}_{oc}(M, [g])$ (see \eqref{eq:defi1} and Definition \ref{defyp}), and Proposition \ref{sameinf}, we have
$$
|||T_h(u)||_{L^{2*}(M, h)}\leq \frac{1}{\sqrt{\mathcal{Y}_{oc}(M, [g])}}||u||_h=\frac{1}{\sqrt{\mathcal{Y}(M, [g])}}||u||_h.
$$
\end{pf}
\vspace{7pt}

\noindent
When $(M, g)=(\mathbb{S}^n, g_{\mathbb{S}^n})$, we have the following obstacle analogue of Lemma \ref{sobolev2}.
\begin{pro}\label{sobolev2op}
Assuming $(M, g)=(\mathbb{S}^n, g_{\mathbb{S}^n})$\; and \;$h=g_w\in [g]$, then for  \;$u\in H^1_+(M, h)$,
\begin{equation}\label{inequality}
||T_h(u)||_{L^{2*}(M, h)}\leq \frac{1}{\sqrt{\mathcal{Y}(\mathbb{S}^n, [g_{\mathbb{S}^n}])}}||u||_h,
\end{equation}
with equality holds if and only if $$u\in C^{\infty}_+(M, h)\;\;\text{and}\;\;R_{g_{uw}}=const.$$ 
\end{pro}

\begin{pf}
The inequality \eqref{inequality} follows directly from Lemma \ref{sobolev1op}. The rigidity part follows from Proposition \ref{sameinf},  Proposition \ref{sameminimizer} and Lemma \ref{sobolev2}. Indeed, $u$\; attains equality in \eqref{inequality} is equivalent to 
$$
||T_h(u)||_{L^{2*}(M, h)}=\frac{1}{\sqrt{\mathcal{Y}(\mathbb{S}^n, [g_{\mathbb{S}^n}])}}||u||_h. 
$$
Thus, Proposition \ref{sameinf} implies \;$u$\; attains equality in \eqref{inequality} is equivalent to 
$$
||T_h(u)||_{L^{2*}(M, h)}=\frac{1}{\sqrt{\mathcal{Y}_{oc}(\mathbb{S}^n, [g_{\mathbb{S}^n}])}}||u||_h.
$$
Now, using the definition of \;$I^{g_{\mathbb{S}^n}}$ (see \eqref{eq:defi1}), we have \;$u$\; attains equality in \eqref{inequality} is equivalent to $$I^{g_{\mathbb{S}^n}}(u)=\mathcal{Y}_{oc}(\mathbb{S}^n, [g_{\mathbb{S}^n}]).$$ 
Moreover, using Proposition \ref{sameminimizer}, we infer that \;$u$\; attains equality in \eqref{inequality} is equivalent $$J^{g_{\mathbb{S}^n}}(u)=\mathcal{Y}(\mathbb{S}^n, [g_{\mathbb{S}^n}]).$$
Thus, using the definition of $J^{g_{\mathbb{S}^n}}$ (see \eqref{eq:defj1}), we have \;$u$\; attains equality in \eqref{inequality} is equivalent  to \;$u$\; attains equality in \eqref{inequality0}. Hence the rigidity part in Lemma \ref{sobolev2} implies \;$u$\; attains equality in \eqref{inequality} is equivalent \;$u\in C^{\infty}_+(M, h)\;\;\text{and}\;\;R_{g_{uw}}=const.$

\end{pf}
\vspace{7pt}

\noindent
\begin{pfn}{ of Theorem \ref{sphere}}\\
It follows directly from Propositon \ref{sameinf}, and the rigidity part in Proposition \ref{sobolev2op}. Indeed, Propositon \ref{sameinf} implies
$I^g(u)=\mathcal{Y}_{oc}(M, [g])$ is equivalent to $I^g(u)=\mathcal{Y}(M, [g])$. Thus, using the definition of $I^g$, we have $I^g(u)=\mathcal{Y}_{oc}(M, [g])$ is equivalent to  $u$ attains equality in \eqref{inequality} with $g=h$. Hence, using the  rigidity part of Proposition \ref{sobolev2op}, we have $I^g(u)=\mathcal{Y}_{oc}(M, [g])$ is equivalent $u\in C^{\infty}_+(M, g)\;\;\text{and}\;\;R_{g_{u}}=const.$
\end{pfn}

\end{document}